\documentclass[12pt]{amsart}
\usepackage{amssymb,epic,eepic,verbatim,graphicx,diagrams}
\input psfig.sty

\newtheorem{theorem}{Theorem}[subsection]

\newtheorem{proposition}[theorem]{Proposition}
\newtheorem{lemma}[theorem]{Lemma}
\newtheorem{lem}[theorem]{}
\theoremstyle{definition}
\newtheorem{definition}[theorem]{Definition}
\theoremstyle{remark}

\newtheorem{example}[theorem]{Example}
\newcommand{\blem}{\begin{lem} \rm}
\newcommand{\elem}{\end{lem}}
\newcommand\mL{\mathcal{L}}

\renewcommand{\L}{\mathcal{L}}

\newcommand{\R}{\mathbb{R}}

\newcommand{\cC}{\mathcal{C}}

\newcommand{\Z}{\mathbb{Z}}

\newcommand\lie[1]{\mathfrak{#1}}

\newcommand{\su}{\lie{su}}
\newcommand{\so}{\lie{so}}

\newcommand{\on}{\operatorname}
\newcommand{\Edge}{\on{Edge}}
\newcommand{\Ver}{\on{Vert}}
\newcommand{\Tet}{\on{Tet}}
\newcommand{\Tri}{\on{Tri}}

\newcommand{\Face}{\on{Face}}

\newcommand{\Hom}{ \on{Hom}}

\newcommand{\SO}{ \on{SO}}

\newcommand{\Vol}{  \on{Vol}}

\newcommand{\ssm}{\kern-.5ex \smallsetminus \kern-.5ex}

\newcommand\dirac{/\kern-1.2ex\partial} 
\newcommand\qu{/\kern-.7ex/} 
\newcommand\lqu{\backslash \kern-.7ex \backslash} 

\newcommand\dr{r_+ \kern-.7ex - \kern-.7ex r_-}
 
\newcommand{\labell}\label

\renewcommand{\d}{{\mbox{d}}}
\newcommand{\ol}{\overline}

\newcommand{\f}{\frac}

\newcommand{\hh}{{\f{1}{2}}}
\newcommand{\thh}{{\f{3}{2}}}
\newcommand{\qq}{{\f{1}{4}}}

\newcommand{\ti}{\tilde}

\newcommand\vol{\on{vol}}

\newcommand\Tr{\on{Tr}}

\newcommand\Map{\on{Map}}

\newcommand\Vect{\on{Vect}}

\newcommand\bdefn{\begin{definition}}
\newcommand\edefn{\end{definition}}
\newcommand\bea{\begin{eqnarray*}}
\newcommand\eea{\end{eqnarray*}}
\newcommand\bcv{\left[ \begin{array}{r} }
\newcommand\ecv{\end{array} \right] }

\newcommand\bma{\left[ \begin{array} }
\newcommand\ema{\end{array} \right]}
\newcommand\ben{\begin{enumerate}}
\newcommand\een{\end{enumerate}}
\newcommand\bex{\begin{example}}
\newcommand\bsj{\left\{ \begin{array}{rrr} }
\newcommand\esj{\end{array} \right\}}

\newcommand\Hess{H}

\newcommand\eex{\end{example}}
\newcommand\sx{*\kern-.5ex_X}

\def\mathunderaccent#1{\let\theaccent#1\mathpalette\putaccentunder}
\def\putaccentunder#1#2{\oalign{$#1#2$\crcr\hidewidth \vbox
to.2ex{\hbox{$#1\theaccent{}$}\vss}\hidewidth}}

\begin{document}

\title{Spherical tetrahedra and invariants of $3$-manifolds}

\author{Yuka U. Taylor}

\address{Department of Mathematics, The George Washington University, 
Washington, DC 20052, 
U.S.A.}\email{yukatylr@gwu.edu}

\author{Christopher T. Woodward}

\address{Mathematics-Hill Center, Rutgers University,
110 Frelinghuysen Road, Piscataway, NJ 08854-8019, U.S.A.}
\email{ctw@math.rutgers.edu}


\maketitle


\section{Introduction}

Let $Y$ be an oriented closed three-manifold and $r$ a positive integer.  The
Reshetikhin-Turaev invariant $Z(Y,r)$ and Turaev-Viro invariant
$TV(Y,r)$ are three-manifold invariants that attempt to make rigorous
the Hamiltonian formulation of quantum Chern-Simons theory.  $Z(Y,r)$
is constructed using the $R$-matrix of the quantized enveloping
algebra $U_q(\lie{sl}_2)$ and Kirby moves, while $TV(Y,r)$ is based on
the $6j$ symbols for $U_q(\lie{sl}_2)$ and a choice of triangulation.
Turaev \cite{vt:90} and Roberts \cite{jr:tv} independently showed that
the $TV(Y,r)$ is the square of the modulus of $Z(Y,r)$.

On the other hand, the Lagrangian (path integral) formulation of
quantum Chern-Simons theory leads to perturbative invariants developed
in \cite{as:cs,as:cs2,le:un}.  The leading term for the invariant is
conjectured in \cite{fg:91} to involve the torsion, the Chern-Simons
invariant, and the spectral flow for flat $SU(2)$ bundles on $Y$.  An
interesting mathematical problem is whether the two formulations can
be shown to agree.  A proof that the leading term is the same for lens
spaces and torus bundles was given in Jeffrey \cite{je:cs}.  Yoshida
\cite{ty:wc} recently announced a proof of equality of the leading
term for a rational homology sphere, using a different definition of
$Z(Y,r)$.

In this paper we apply our previous work on asymptotics of the quantum
$6j$ symbols \cite{tw:6j} to the asymptotics of $TV(Y,r)$ as $r \to
\infty$.  Substituting the asymptotic formula and applying stationary
phase yield a finite dimensional integral involving Gram matrices of
spherical tetrahedra which turns out to be a spherical version of an
integral considered by Ponzano-Regge \cite{pr:68} and Korepanov
\cite{ik:3d}; see also Mizoguchi and Tada \cite{mi:th}.
Unfortunately, we have nothing rigorous to say about the asymptotics
because of various problems involving convergence of the integral and
error estimates for the asymptotics of the $6j$ symbols.  The modest
results of this paper are a proof that the integral is invariant under
the Pachner moves, as one would expect from the connection with
Turaev-Viro, and of convergence for the sphere $S^3$.

We would like to thank I. Korepanov, M. Leingang, F. Luo, I. Rivin,
and J. Roberts for their comments and discussions.  This material
appeared in the first author's Rutgers University 2003 Ph.D. thesis.

\section{$6j$ Symbol for $U_q(\lie{sl}_2)$ and its asymptotic formula}\label{intro6j}

Let $U_q(\lie{sl}_2)$ denote the quantized enveloping algebra at a primitive
$2r$-th root of unity $q = \exp(\pi i /r)$.  Let $[n]_q$ be the {\it
quantum integer $n$} defined by 
$$[n]_q=\frac{q^n-q^{-n}}{q-q^{-1}} = \frac{ \sin( n \pi / r)}{\sin (
  \pi/r)} $$ 
for $n \in \mathbb{Z}$.  We say that a half-integer $j$ 
is a {\em color} at level $r$ if 
$$ 0 \leq j \leq \frac{r-2}{2} .$$
For any color $j$, define
$$ \Delta_j = (-1)^{2j} [2j + 1].$$
A triple of colors $j_1,j_2,j_3 \in \Z/2$ is called {\em admissible}
if
$$ \max(j_1 - j_2,j_2 - j_1) \leq j_3 \leq \min(j_1 + j_2,r-2 - j_1 -
j_2) $$
and
$$ j_1 + j_2 + j_3 \in \Z .$$
The quantity 
$$ \Delta = \Delta_a^{-1} \sum_{b,c, (a,b,c) \ \text{admissible}}
\Delta_b \Delta_c = r \sin(\pi/r)^{-2} $$
\cite{tv:6j} and in particular is independent of $a$.

For any $6$-tuple of colors $j_{ab}, 1 \leq a < b \leq 4$, the {\em
quantum $6j$ symbol}
$$ \bsj j_{12} & j_{23} & j_{13} \\
        j_{34} & j_{14} & j_{24} \esj $$
is a rational number obtained from associativity of the tensor product for
representations of $U_q(\lie{sl}_2)$.  There are two standard conventions
for $U_q(\lie{sl}_2)$ with tetrahedral symmetry, which are related by a
sign
$$ (-1)^{\sum_{a<b} 2j_{ab}} .$$
The reader should note that the Turaev-Viro convention is different
from the convention in our earlier paper \cite{tw:6j}, which was
chosen because it agrees with the accepted conventions for $q = 1$.
The $6j$ symbol satisfies the {\em orthogonality relations}
\cite{cq:95}
\begin{equation}\label{normal}
 \sum_{j_{14}} \Delta_{j_{14}} \Delta_{m}
\bsj j_{12} & j_{13} & n \\
                j_{34} & j_{24} & j_{14} \esj
\bsj j_{12} & j_{13} & m \\
                j_{34} & j_{24} & j_{14} \esj=\delta_{m,n}, 
\end{equation}
and the {\em pentagon} or Biedenharn-Elliot relation
\begin{equation}\label{pentagon}
\tau(1234) \} \{ \tau(2345) \} = \sum_{j_{15}}(-1)^z[2j_{15}+1] \{
\tau(1235) \} \{ \tau(1345) \} \{ \tau(1245) \}
\end{equation}
where $z$ is the sum of all $j_{ab}$, $a,b \in \{ 1,2,3,4,5 \}$ and
$(j_{23},j_{34},j_{24})$ is $q$-admissible, and 
$\{ \tau(abcd) \}$ is short for 
$$ \{ \tau(abcd) \} = 
\bsj j_{ab} & j_{bc} & j_{ac} \\
     j_{cd} & j_{ad} & j_{bd} \esj .$$

In our previous paper \cite{tw:6j} we obtained the following result on
the asymptotics of the quantum $6j$ symbols as the labels and level
are simultaneously rescaled.  Set
$$r(k) \equiv k(r-2)+2.$$ 
Let $\tau$ denote the tetrahedron in the sphere $S^3$ with edge
lengths
\begin{equation}\label{l}
l_{ab}=2\pi\left(\frac{kj_{ab} + \hh}{r(k)}\right),
\end{equation}
if it exists, and let $\theta_{ab}$ denote the exterior dihedral
angles.  Define 
$$\phi= \frac{r(k)}{2\pi} \left(\sum_{a < b}l_{ab}\theta_{ab} -
2\vol(\tau) \right) $$
and
$$ G(\tau) = \det(\cos(l_{ab})) $$ 
where $(\cos(l_{ab}))$ is the spherical $4 \times 4$ Gram matrix.
Then,
\begin{equation}\label{S^3-asymp}
\bsj kj_{12} & kj_{13} & kj_{23} \\
                kj_{34} & kj_{24} & kj_{14} \esj_{q} \sim
\frac{2\pi \cos(\frac{\pi}{4}+\phi)}
{(r(k))^{\frac{3}{2}}G^{\frac{1}{4}}(\tau(l_{ab}))},
\end{equation} 
if $\tau$ exists and is non-degenerate.

\section{The Turaev-Viro invariant}

Let $Y$ be a compact triangulated 3-manifold with tetrahedra
$\Tet(Y)$, triangles $\Tri(Y)$, edges $\Edge(Y)$, and vertices
$\Ver(Y)$.  A {\em coloring} of $Y$ at an integer $r \ge 2$ is a map
$$j:\Edge(Y) \to \left\{ 0, \hh, \ldots,
\frac{r-2}{2}\right\} .$$ 
For each such coloring, define
$$ TV(Y,r,j) = \Delta^{- v(Y)} \prod_{e \in \Edge(Y)} \Delta_{j(e)}
\prod_{\tau \in \Tet(Y)} \{ j(\tau) \}_q $$
$j(\tau)$ denotes the vector of values of $j$ on the 6 edges of a
tetrahedron $\tau$, $\{ j(\tau) \}_q$ is the $6j$-symbol for
$U_q(\lie{sl}_2), q = \exp ( \pi i /r)$ for the colors associated to
the edges of the tetrahedron $\tau$.  The Turaev-Viro invariant of $Y$
is
\begin{equation} \label{tv} TV(Y,r) = \sum_{j} TV(Y,r,j).
\end{equation}
The pentagon and orthogonality identities imply that $TV(Y,r)$ is
invariant under the Pachner 2-3 and 1-4 moves and hence
independent of the triangulation, that is, a topological invariant of
$Y$.  

\section{Non-Euclidean Tetrahedra}\label{geomfact}

This section provides various elementary facts about non-Euclidean
tetrahedra relevant to this paper. For the proof, we refer to
\cite{tw:6j}.  Let $E^n, S^n$ denote $n$-dimensional Euclidean,
spherical space respectively.  Let $S_n$ denote an $n$-dimensional
simplex and $l_{ab}$ edge lengths in $S_n$.  The {\it Cayley-Menger
determinant} for a Euclidean simplex $S_n$, denoted by $G_0(l_{ab})$, is
defined by
\begin{equation}\label{cayley}
G_0(l_{ab})=\det
\left( \begin{array}{rrrrrr}
0 & 1  & 1 & 1 &\ldots & 1 \\
1 & 0  & -\hh l_{12}^2 & -\hh l_{13}^2 & \ldots  &-\hh l_{1n}^2 \\
1 &-\hh l_{21}^2 &  0   &-\hh l_{23}^2 & \ldots&-\hh {l_{2n}}^2 \\
\vdots & \vdots & \vdots &  & \ldots & \vdots\\
1 & -\hh l_{n1}^2 & -\hh l_{n2}^2 & \ldots & \ldots  & 0
\end{array}\right).
\end{equation}
For a spherical simplex, we define $n \times n$ {\it Gram matrix}
\[ G(l_{ab})=\det (\cos(l_{ab}))  .\]
Note that this is the volume of a Euclidean $n+1$-simplex with $n$
vertices on the unit sphere and one at $0$.  We will need later the
following facts on Euclidean and spherical simplices; the hyperbolic
versions are discussed in \cite{tw:6j}.

\begin{theorem}\label{misc}

\begin{enumerate}

\item (Cayley formula, \cite[p. 98]{Bl:70}) If a simplex $S$
with edge lengths $l_{ab}$ exists in $E^n$, then
\[ (n!\Vol(S))^2= G_0(l_{ab}). \]
\item (Schl\"{a}fli formula, \cite[p. 281]{jm:94})
\label{shlafli} For
an $n$-simplex $S$ in $E^n$ or $S^n$,
\[ (n-1)\kappa \d \Vol_n(S)=\sum \Vol_{n-2}(F) \d \theta_F \]
qwhere the sum is over $(n-2)$-dimensional faces $F$ of the simplex
$S$, $\theta_F$ is the exterior dihedral angle around $F$ and $\kappa
= 0 , 1$ is the curvature.

\item
For the case of a triangle, one has the factorizations:
\begin{align*}
G_0 &=\frac{1}{4}(l_{12}+l_{23}+l_{13})(l_{12}+l_{23}-l_{13})
(l_{12}-l_{23}+l_{13})(-l_{12}+l_{23}+l_{13})\\
G &=4\sin(\hh(l_{12}+l_{23}+l_{13}))\sin(\hh
(l_{12}+l_{23}-l_{13}))\sin(\hh
(l_{12}-l_{23}+l_{13}))\\
&\sin(\hh(-l_{12}+l_{23}+l_{13}))
\end{align*}

\item
A Euclidean triangle exists if and only if
\begin{equation}\label{triangle}
l_{12} \leq l_{13}+l_{23}, \ l_{13} \leq l_{12} + l_{23} \ 
l_{23} \leq l_{12} + l_{13}. 
\end{equation}
A spherical triangle exists if and only if \eqref{triangle} and
\[ l_{12}+l_{13}+l_{23} \leq 2 \pi .\]  
\item \label{tetra} A non-degenerate tetrahedron with edge lengths
$l_{ab}$ exists in $E^3$, $S^3$ respectively if and only if $l_{ab}$
satisfy \eqref{triangle} for faces and $ I > 0, \ G > 0 $
respectively.
\item \label{wigner} The derivative of an edge length $l_{ab}$ in a Euclidean
resp. spherical tetrahedron $\tau$ with respect to an opposite
dihedral angle $\theta_{cd}$ is given by
\[ 
\frac{\partial l_{ab}}{\partial \theta_{cd}}
= \pm \frac{G_0^{1/2}(l_{ij})}{l_{ab} l_{cd}}, \ \ \ \ 
\frac{\partial l_{ab}}{\partial \theta_{cd}}
= \pm \frac{G^{1/2}(l_{ij})}{\sin(l_{ab})\sin(l_{cd})},
\]
\end{enumerate}
\end{theorem}

\section{Asymptotic Pentagon and Normalization Identities}\label{apo}

In this section, we prove several geometric identities which may be
viewed as semiclassical analogs of the identities \eqref{pentagon},
\eqref{normalization} for $6j$ symbols. We will use them when we
discuss $3$-manifold in section \ref{appl}.  The Euclidean versions
are due to Ponzano and Regge \cite{pr:68}.  Starting from this
section, we fix an integer $r \geq 3$ and $q=\exp(\frac{\pi i}{r})$.

A simplex spanned by vertices $v_0, \ldots, v_n$ is denoted by $S_{0
\ldots n}$.  If the simplex is two-dimensional, then we sometime
denote $S_{kl}$ by $e_{kl}$.  Also, a vector from  $v_i$ to 
$v_j$
is denoted by $v_{ij}$.  A simplex spanned by vertices $v_0,
\ldots, \hat{v_i}, \ldots, v_n$, in which the vertex $v_i$ is omitted,
is denoted by $S_{i}$.  A volume of a simplex is written as 
$\vol(S_{0 \ldots n})$, or $\vol(S_{h_1,\ldots, h_n})$,where $h_i$ are
spanning vectors of the simplex.  
If a simplex is two-dimensional, we sometimes write $l_{kl}$ for 
the length $\vol(S_{kl})=\vol(e_{kl})$.

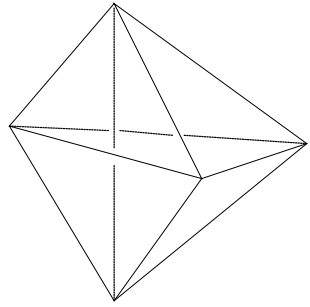
\begin{figure}[h]
\begin{center}
\setlength{\unitlength}{0.0003in}
\begingroup\makeatletter\ifx\SetFigFont\undefined%
\gdef\SetFigFont#1#2#3#4#5{%
  \reset@font\fontsize{#1}{#2pt}%
  \fontfamily{#3}\fontseries{#4}\fontshape{#5}%
  \selectfont}%
\fi\endgroup%
{\renewcommand{\dashlinestretch}{30}
\begin{picture}(5124,5139)(0,-10)
\path(1812,5112)(12,3012)(3312,2112)
        (1812,5112)(5112,2712)(3312,2112)
        (1812,12)(12,3012)
\path(5112,2712)(1812,12)
\dottedline{45}(1812,5037)(1812,2637)
\dottedline{45}(1812,2337)(1812,12)
\dottedline{45}(12,3012)(1712,2937)
\dottedline{45}(1912,2937)(2812,2850)
\dottedline{45}(3012,2850)(5112,2712) 
\end{picture}
}

\end{center}
\caption{2-3 move}\label{23}
\end{figure}

Consider a complex with vertices $v_0, \ldots, v_4$ and five
tetrahedra $S_{0}$, $S_{4}$, $S_{1}$, $S_{2}$, $S_{3}.$ Suppose the
complex described above is embedded in $S^3$.  For the spherical Gram
matrix associated with $S_{i}$, let $G_i$ denote its determinant.  For
$j = 1,\ldots,5$ define $s_j = 1$ resp. $-1$ if the embedding is
orientation preserving resp. reversing.  Around the edge $e_{04}$, we
have three exterior dihedral angles $\theta^1_{04}$, $\theta^2_{04}$,
$\theta^3_{04}$.  Define the {\em defect angle} around the edge
$e_{04}$ by
\[ \omega_{04}= \sum_{j=1}^3  s_j (\pi - \theta^j_{04})  .\]

\begin{theorem}  The determinants $G_i$ of spherical Gram matrices for
tetrahedra $S_i$ satisfy the identities
\begin{enumerate}
\item(Asymptotic pentagon identity) In the same situation as in
\eqref{pentagon}
\begin{equation}\label{sjac}
\frac{\partial \omega_{04}}{\partial l_{04}}= s_1 s_2 s_3
\sin^2(l_{04}) \sqrt{\frac{ G_0 G_4}{ G_1G_2G_3}},
\end{equation} 
\item(Asymptotic normalization identity) 
In the same situation described as in
\eqref{normal}
\begin{equation}\label{normalization}
\sin (l_{cd}) \int \frac{\sin (l_{ab})}{\sqrt{G(\tau(l_{ij}))}}
\ d l_{ab}= \pi. 
\end{equation}
Here, $l_{ab}$ and $l_{cd}$ are the lengths of opposite edges and
$G(\tau(l_{ab}))$ is the determinant of the spherical Gram matrix
associated with the tetrahedron with edge lengths $l_{ij}$.
\end{enumerate}
\end{theorem}

The equation \eqref{sjac} can be obtained heuristically via stationary
phase applied to \eqref{pentagon}.  To prove \eqref{normalization},
note that
\[ \int \frac{\partial \theta_{cd}}{\partial l_{ab}}\d
l_{ab}=\int_0^\pi \d \theta_{cd}=\pi. \]
by Theorem \ref{misc} \eqref{wigner}.  The proof of \eqref{sjac} is by a
series of lemmas.  Suppose $v_0, \ldots, v_n$ in $\mathbb{R}^n$ form a
$n-$simplex $S_{0,\ldots,n}$.  Consider the $(n-3)-$simplex
$S_{0,4,5,\ldots, n}$.  Let $h_i$ be the vector starting at the vertex
$v_i$ perpendicular to the simplex $S_{0,4,5,\ldots,n}$ for each $i =
1,2,3.$

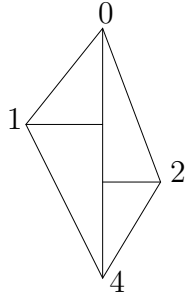
\begin{figure}[h]
\begin{center}
\setlength{\unitlength}{0.00033333in}
\begingroup\makeatletter\ifx\SetFigFont\undefined%
\gdef\SetFigFont#1#2#3#4#5{%
  \reset@font\fontsize{#1}{#2pt}%
  \fontfamily{#3}\fontseries{#4}\fontshape{#5}%
  \selectfont}%
\fi\endgroup%
{\renewcommand{\dashlinestretch}{30}
\begin{picture}(2642,4350)(0,-10)
\path(1500,4125)(300,2625)(1500,225)
        (2400,1725)(1500,4125)
\put(1425,4200){\makebox(0,0)[lb]{\smash
{{{\SetFigFont{12}{14.4}{\rmdefault}{\mddefault}{\updefault}0}}}}}
\put(0,2550){\makebox(0,0)[lb]{\smash
{{{\SetFigFont{12}{14.4}{\rmdefault}{\mddefault}{\updefault}1}}}}}
\put(1600,0){\makebox(0,0)[lb]{\smash
{{{\SetFigFont{12}{14.4}{\rmdefault}{\mddefault}{\updefault}4}}}}}
\put(2550,1725){\makebox(0,0)[lb]{\smash
{{{\SetFigFont{12}{14.4}{\rmdefault}{\mddefault}{\updefault}2}}}}}
\path(1500,4125)(1500,225)(1500,4125)
\path(300,2625)(1500,2625)(300,2625)
\path(1500,1725)(2400,1725)(1500,1725)
\end{picture}
}
\end{center}
\caption{The vectors $h_1$ and $h_2$}\label{h1h2}
\end{figure}
The length $\Vert h_i \Vert$ of the vector is the distance from $v_i$
to $S_{0,4,5,\ldots,n}$.  Also, the dihedral angle between
$3-$dimensional simplices $S_{0,i,j,4}$ and $S_{0,j,k,4}$ for $i \neq
j \neq k \in \{ 1,2,3 \}$ around $S_{0,4,\ldots,n}$ is the same as the
angle $\phi_{ik}$ between vectors $h_i$ and $h_k$.  In particular, the
exterior dihedral angle between $S_{0,i,j,4}$ and $S_{0,j,k,4}$ is
$\pi - \phi_{ik}$, which we denote by $\theta_{ik}$.  The volume of
the $(n-1)$-simplex $S_i=(v_0,\ldots, \hat{v_i}, \ldots,v_n)$, denoted
by $V_i$, is
\begin{eqnarray}\label{volume}
V_i&=&\frac{1}{(n-1)(n-2)} \Vol(0,h_j,h_k) \Vol(S_{0,4,\ldots,n}) \\
&=& \frac{1}{(n-1)(n-2)}\Vert h_j\Vert \Vert h_k \Vert
\Vol(S_{0,4,\ldots,n})\sin \theta_{jk}
\end{eqnarray}
where $i,j,k \in \{ 1,2,3 \}, i \neq j \neq k$.  The volume of the
$n$-simplex $S_{0,\ldots,n}$ is
\begin{equation}\label{H and I}
\binom{n}{3}\Vol(S_{0,\ldots,n})=\Vol(h_1,h_2,h_3)\Vol(S_{0,4,\ldots,n}).
\end{equation}
Indeed, 
\begin{align*}
n!\Vol(S(0,\ldots,n))& =\det (e_{01}, e_{02}, e_{03},e_{04}, \ldots,
e_{0n})\\ &=\det (h_1,h_2, h_3)\det(e_{04},\ldots, e_{0n})\\ &=
3!(n-3)!\Vol(h_1,h_2,h_3)\Vol(S_{0,4,\ldots,n}).
\end{align*}
\begin{lemma}
\begin{equation} \label{H}
\frac{\partial (\Vol(h_1,h_2,h_3))^2}{\partial \omega_{04}}
|_{\omega_{04} = 0} =\frac{2(n-1)^3(n-2)^3}{(3!)^2}\frac{s_1 s_2 s_3
  V_1V_2V_3}{(\Vol(S_{0,4,\ldots,n})) ^3}
\end{equation}
\end{lemma} 

\begin{proof} We know that
$$\Vol(h_1,h_2,h_3) = \frac {1}{3!} \det(h_i \cdot h_j)^{1/2} = \frac
{1}{3!}\Vert h_1 \Vert \Vert h_2\Vert \Vert h_3 \Vert {\det
(\cos\phi_{ij})}^\frac{1} {2} .$$
Substituting $\phi_{ij}=\pi - \theta_{ij}$ and expanding the
determinant yield
$$(\Vol(h_1,h_2,h_3))^2=(\frac{1}{3!})^2(\Vert h_1 \Vert \Vert h_2 \Vert \Vert h_3 \Vert)^2(\cos^2
\theta_{12}+\cos^2 \theta_{13}+ \cos^2 \theta_{23}+2\cos \theta_{12}
\cos \theta_{13} \cos \theta_{23}).$$

The differential with respect to $\theta_{12}, \theta_{13},
\theta_{23}$ is

$$\d (\Vol(h_1,h_2,h_3))^2=(\frac{1}{3!})^2(\Vert h_1 \Vert \Vert h_2 \Vert \Vert h_3 \Vert)^2 \{
2\sin \theta_{12}(\cos
\theta_{12}+ \cos \theta_{13}\cos
\theta_{23}) \d \theta_{12}+$$
$$2\sin \theta_{13}(\cos
\theta_{13}+ \cos \theta_{12}\cos
\theta_{23}) \d \theta_{13}+$$
$$2\sin \theta_{23}(\cos
\theta_{23}+ \cos \theta_{12}\cos
\theta_{13}) \d \theta_{23} \}.$$

The double angle formula, together with $\omega_{04} = 0$ gives

$$\cos \theta_{12}=\cos (\pi - (s_{2} \theta_{13}+ s_{1}
\theta_{23}))=-\cos \theta_{13}\cos \theta_{23}+ \sin
s_{2} \theta_{13}\sin s_{1} \theta_{23}.$$

Therefore, $\d(\Vol(h_1,h_2,h_3))^2\mid_{\omega_{04} = 0}$ is equal 
to 
$$ \frac {2}{(3!)^2}
(\Vert h_1 \Vert \Vert h_2 \Vert \Vert h_3 \Vert)^2\sin s_{3}
\theta_{12} \sin s_{2} \theta_{13}\sin s_{1} \theta_{23}(s_{3}
d\theta_{12}+s_{2} d\theta_{13}+s_{1} d\theta_{23}).$$
By  (\ref{volume}),
$$V_1V_2V_3=(\frac {1}{(n-1)(n-2)})^{3}(\Vert h_1 \Vert \Vert h_2 \Vert \Vert h_3 \Vert)^2(\Vol(S_{0,4,\ldots,n}))^3
\sin \theta_{12}\sin  \theta_{23}\sin  \theta_{13}.$$
The lemma follows since 
$\d\omega_{04} = \sum_{k \neq i \neq j} s_{k} \d\theta_{ij} .$
\end{proof}

Let $x$ be the length of the edge $e_{ij}$ from the vertex $v_i$ to
the vertex $v_j$ and $x_{\pm}$ be the roots of the Cayley-Menger 
determinant associated with the $n$-simplex.  
\begin{lemma}
\begin{equation}\label{dI over dx}
\frac {\partial \Vol(S_{0,\ldots,n})^2}{\partial
x^2}\vert_{x^2=x^2_{\pm}}={\pm} \frac{1}{n^2}V_iV_j,
\end{equation}

\end{lemma}

\begin{proof}  
Without loss of generality, assume $x$ is the length of $e_{0n}$.
Using the Cayley-Menger determinant,
$$ \Vol^2(S_{0,\ldots,n}) = \frac{-1}{4(n(n-1))^2}
\Vol^2(S_{1,\ldots,n-1})
(x^2- x_-^2)(x^2 - x_+^2) $$
and so 
\begin{equation} \label{part}
 \frac{\partial \Vol^2(S_{0,\ldots,n})}{\partial x^2}
\vert_{x^2=x^2_{\pm}}=
\frac{\pm 1}{4(n(n-1))^2} \Vol^2(S_{1,\ldots,n-1}) (x_+^2 - x_-^2).
\end{equation}
The roots $x_\pm$ correspond to values of the length for which the
simplex embeds into $\R^n$.  We can choose the embeddings so that only
the image $v_{n+},v_{n-}$ of the vertex $v_n$ varies; see Figure
\ref{n=3}.  Let $h_0$, resp. $h_n$, denote distance of $v_0$,
resp. $v_{n,\pm}$, to $S(1,\ldots,n-1)$, so that
$$ \Vol(S_{n}) = \frac{1}{n-1} \Vol(S_{1,\ldots,n-1}) h_0, \ \
\ \ \Vol(S_{0}) = \frac{1}{n-1} \Vol(S_{1,\ldots,n-1}) h_n.
 $$
\begin{figure}
\begin{center}
\setlength{\unitlength}{0.00055000in}
\begingroup\makeatletter\ifx\SetFigFont\undefined%
\gdef\SetFigFont#1#2#3#4#5{%
  \reset@font\fontsize{#1}{#2pt}%
  \fontfamily{#3}\fontseries{#4}\fontshape{#5}%
  \selectfont}%
\fi\endgroup%
{\renewcommand{\dashlinestretch}{30}
\begin{picture}(4217,3750)(0,-10)
\path(2100,3375)(300,1275)(3900,1275)(2100,3375)
\path(300,1275)(2700,375)(3900,1275)
\path(300,1275)(2700,2175)(3900,1275)
\dottedline{45}(2100,3375)(2100,1275)
\dottedline{45}(2700,375)(2700,1275)
\dottedline{45}(2100,3375)(2700,2175)
\dottedline{45}(2100,3375)(2700,375)
\put(2100,3600){\makebox(0,0)[lb]{\smash{{{\SetFigFont{11}{13.2}
{\rmdefault}{\mddefault}{\updefault}0}}}}}
\put(0,1200){\makebox(0,0)[lb]{\smash{{{\SetFigFont{11}{13.2}
{\rmdefault}{\mddefault}{\updefault}1}}}}}
\put(4125,1200){\makebox(0,0)[lb]{\smash{{{\SetFigFont{11}{13.2}
{\rmdefault}{\mddefault}{\updefault}2}}}}}
\put(2700,0){\makebox(0,0)[lb]{\smash{{{\SetFigFont{11}{13.2}
{\rmdefault}{\mddefault}{\updefault}3+}}}}}
\put(2625,1875){\makebox(0,0)[lb]{\smash{{{\SetFigFont{11}{13.2}
{\rmdefault}{\mddefault}{\updefault}3-}}}}}
\end{picture}
}
\end{center}
\caption{n=3}\label{n=3}
\end{figure}

Let $w$ denote the distance from the projection of $v_n$, to the 
projection of $v_0$ in $S(1,\ldots,n-1)$. By the Pythagorean theorem,
$$ x_+^2 = (h_0 + h_n)^2 + w^2, \ \ x_-^2 = (h_0 - h_n)^2 + w^2 .$$
Hence $x_+^2 - x_-^2 = 4 h_0 h_n$, so the lemma follows from
\eqref{part}.
\end{proof}

Finally we prove the asymptotic pentagon identity.  We use the above lemmas for $n=5$.
Suppose that the vertices $v_0, v_1, v_2, v_3, v_4$ lie in $S^3$ 
and $v_5=0$.  Let
$$ I = \Vol^2(S_{0,\ldots,4}).$$

It suffices to compute
$$\frac{\partial \omega_{04}}{\partial y}= \frac{\partial \omega_{04}}{\partial
I}\frac{\partial I}{\partial x^2}\frac{\partial x^2}{\partial y},$$
where $x$ is the Euclidean length between $v_0$ and $v_4$ and $y$ the
spherical geodesic distance.  By \eqref{H and I},
$$I =\frac{\Vol^2(S_{0,4,5})\Vol^2(h_1,h_2,h_3)}{5^22^2}.$$
Because $\Vol(S_{0,4,5})$ is independent of $\omega_{04}$,
$$\frac{\partial I}{\partial
\omega_{04}}=\frac{\Vol^2(S_{0,4,5})}{5^22^2}\frac{\partial
(\Vol^2(h_1,h_2,h_3))}{\partial \omega_{04}}.$$
By \eqref{H},
$$ \frac{\partial \omega_{04}}{\partial
I}=\frac{5^22^2}{96}\frac{\Vol(S_{0,4,5})s_1s_2 s_3}{V_1V_2V_3}.$$
By \eqref{dI over dx},
$$s_1s_2 s_3\frac{\partial \omega_{04}}{\partial
x^2}=\frac{5^22^2}{96}\frac{\Vol(S_{0,4,5})}{V_1V_2V_3}
\frac{V_0V_4}{5^2}=\frac{1}{24}\Vol(S_{0,4,5})
\frac{V_0V_4}{V_1V_2V_3}.$$
Note that $\Vol(S_{0,4,5})=\frac{1}{2}\sin (y)$ and $x=2\sin
(\frac{y}{2}),$ where $x$ is the length of the straight line
from $v_0$ and $v_4$ and $v_5=0$.  Hence,
$$\frac{\d x^2}{\d y}=4\sin (\frac{y}{2})\cos (\frac{y}{2})= 
2 \sin (y).$$
Thus, 
$$\frac{\partial \omega_{04}}{\partial y}=
\frac{\partial \omega_{04}}{\partial x^2}\frac{\partial x^2}{\partial y}
=\frac{s_1s_2s_3}{24}\frac{\sin(y)}{2}2\sin(y)
\frac{V_0V_4}{V_1
V_2V_3} = 
s_1s_2s_3 \sin^2(y)\sqrt{\frac{G_0G_4}{G_1G_2G_3}}.$$

\section{A semiclassical three-manifold Invariant}\label{appl}

In this section, we explain how to use
\eqref{sjac},\eqref{normalization} to define a formal three-manifold
invariant which is a spherical version of the formal invariant
introduced by Korepanov in \cite{ik:3d} and \cite{ik:ls}.  By formal
we mean that the existence of the invariant depends on the convergence
of certain finite dimensional integrals, which we can only prove in
the case of $S^3$.

\subsection{Definition of the Invariant}\label{def}

Let $Y$ be a triangulated, closed, and oriented three-manifold with
vertices $\Ver(Y)$, edges $\Edge(Y)$, triangles $\Tri(Y)$, and
tetrahedra $\Tet(Y)$.  Let $\mL$ denote the space of the
edge-labellings
$$ \mL = \{ l : \ \Edge(Y) \to [0,\pi], \ \ G(l(\tau)) > 0 \  
\forall \tau \in \Tet(Y) \} .$$
Here, $l(\tau)$ denotes the $6$-tuple which is a restriction of a
labelling $l$ on the edges in $\tau$ and $G(l(\tau))$ the determinant
of the spherical $4 \times 4$ Gram matrix associated with 
$l(\tau)$, and the edge length $l_{ab}$ is as defined in \eqref{l}.  

By Theorem \ref{misc}, if $G(l_{ab})>0$, there is a non-degenerate
spherical tetrahedron with edge length $l_{ab}$.  So, given an $l \in
\mL$ and $\tau \in \Tet(Y)$, there is an embedding $\varphi: \tau \to
S^3$ such that for any edge $e \subset \tau$, the length of the edges
of the tetrahedron $\varphi(e)$ is $l(e)$.  For any coloring $l$ and
any edge $l(e):= l_e$ in the spherical tetrahedron $\varphi(\tau)$,
let $\phi_{l_e,t}$ resp. $\theta_{l_e,t}$ denote the interior
resp. exterior dihedral angle at $l_e$ in $\varphi(\tau)$.  Let 
$$ s: \Tet(Y) \to \{ \pm 1 \} $$
be a sign assignment to each tetrahedron in $Y$.  For each $e \in
\Edge(Y)$ and labelling $l$, define the {\it defect angle around
the edge $e$} to be
\begin{equation}\label{defect}
 \omega_{l_e,s} = 2\pi - \sum_{\tau \supset e} 
s(\tau) \phi_{l_e,\tau} .
\end{equation}
We say that a labelling $l$ is {\em flat} with respect to the sign
choice $s$ if
$$\omega_{l_e,s} = 0 \ \text{mod} \ 2\pi  \ \ \forall e \in \Edge(Y).$$

\begin{definition}
$\mL_{\flat,s}$ denotes the set of
flat labellings  with a fixed sign assignment $s$.  That is, 
$ \mL_{\flat,s}=\{ l \in \L : \omega_{l_e,s} = 0 \ \text{mod} \ 2\pi\}. $
\end{definition}

\begin{proposition}\label{varphi}  
Suppose that $Y$ is simply connected.  For a given flat labelling $l$
and a fixed sign assignment there exists a map $\varphi: \ Y \to S^3$
such that $\varphi|_{\tau}$ is an embedding of $\tau$ with length
$l_{\tau}$, for all tetrahedra $\tau \in \Tet(Y)$.  Any other map
$\varphi': \ Y \to S^3$ whose restriction to a tetrahedron is an
embedding is obtained by composing $\varphi: \ Y \to S^3$ with an
element of $SO(4)$.
\end{proposition}
\noindent The proof is similar to the construction of developing maps for
hyperbolic or spherical manifolds and is left to the reader.

Suppose that $Y$ is not necessarily simply connected.  Let $\ti{Y} \to
Y$ be the universal cover of $Y$.  Each flat labelling $l$ with a
fixed sign assignment $s$ defines $\varphi_l: \ \ti{Y} \to S^3$.  Let
$|\tau |$ denote the spherical tetrahedron $\varphi_l(\tau)$ realized
from $l(\tau)$.  For any $\gamma \in \pi_1(Y)$, $\gamma |\tau| $ is a
spherical tetrahedron, related to $| \tau |$ by an element $\rho
(\gamma)$ in the isometry group $SO(4)$ of $S^3$.  By construction
$$\varphi_l(\gamma_1 \gamma_2 |\tau|) = \rho (\gamma_1) \rho(\gamma_2)
\varphi_l(|\tau|) .$$ 
It follows that $\rho$ is a homomorphism
$$ \rho : \ \pi_1(Y) \to SO(4) = (SU(2) \times SU(2))/ \{ \pm 1 \} .$$
Let $[\rho ]$ denote the conjugacy class of $\rho$ in the representation
variety 
$$ R(Y,SO(4)) := \Hom(\pi_1(Y),SO(4))/SO(4) . $$
Because of the last statement in proposition \ref{varphi}, 
$[\rho]$ is independent of the choice of the base tetrahedron $\tau$ or
an embedding $\tau \to S^3$.  
Let $\mL_{\flat, [\rho]} = \cup_s \mL_{\flat, 
[\rho],s}$ denote
the set of flat labellings $l$ which give rise to the class $[\rho]$.  
 
Given $l \in \mL_{\flat,[\rho],s}$, recall that the defect angle
$\omega_{l_e}$ around
an edge $e$ is defined by \eqref{defect}.    
Let $H$ denote the matrix
$$H = ( \frac{\d \omega_i}{\d l_j} )_{i,j \in \Edge(Y)} .$$
By Schl\"{a}fli's formula \ref{shlafli},
$H$ is the Hessian of the function
$$\sum_{e \in \Edge(Y)} \omega_{l_e,s} l_e - \sum_{\tau \in \Tet(Y)}
s(\tau) 2\vol(|\tau|) ;$$
in particular, $H$ is symmetric.  

For any matrix $M = (m_{ij}), i,j \in \Edge(Y)$, and subsets $I,J
\subset \Edge(Y)$, we denote by $M_{IJ}$ the sub-matrix of $M$
obtained by restricting the index set for rows, resp. columns, to $I$,
resp. to $J$.  Let $\cC \subset \Edge(Y)$ be a maximal subset of edges
such that the sub-matrix $H_{\cC \cC} \subset A$ is positive definite.
Let $\ol{\cC}$ denote its complement $\Edge(Y) \ssm \cC$.  Define
\begin{equation} \label{tvsc}  
I(Y,[\rho]) := (\frac{1}{2\pi})^{\# \Ver} 
\sum_s \int_{l \in
\mL_{\flat,[\rho],s}} \prod_{\tau \in
\Tet}G(l(\tau))^{-1/4}\prod_{e \in \Edge} \sin
(l_e) \frac{\bigwedge_{e \in \ol{\cC}} \d l_e}{\sqrt{\det ( 
H_{\cC \cC})}}.
\end{equation}
If $R(Y,SO(4))$ is finite, then we define
\[ I(Y):=\sum_{[\rho] \in R(Y,SO(4))} I(Y,[\rho]) .\]
This is not exactly the expression predicted by stationary phase
applied to $TV(Y,r)$; that expression is (even) more complicated due
to the inclusion of phases and certain powers of $2$ which we have
ignored.  These omissions are partly discussed in the last section of
the paper.

\subsection{Formal topological invariance}

By Pachner's theorem \cite{pa:tr}, any two triangulations of a given
$3$-manifold are related by a sequence of $1$-$4$ and $2$-$3$ moves.
The $1$-$4$ move replaces a tetrahedron with four tetrahedra by adding
a vertex or vice versa. The Pachner $2$-$3$ move replaces two
tetrahedra sharing a face with three tetrahedra by adding an edge or
vice versa.

\begin{theorem}
$I(Y,[\rho])$ is a formal topological invariant, i.e.,
independent of the choice of $\cC$ and invariant under the Pachner
moves assuming convergence.
\end{theorem}

\noindent 
First we show invariance of the integral under a $2$-$3$ move.  In the
triangulation of $Y$, find a complex of two tetrahedra with vertices
$v_0,v_1,v_2,v_3,v_4$.  Denote it by $X$.  For $\Edge(Y)$, we have the
set of labellings $\mL_\flat=\bigcup_s \mL_{\flat,s},$ where $s$ is a
sign assignment $\Tet(Y) \to \{ \pm 1 \}$.  Consider a new
triangulation $T'$ of $Y$, obtained by adding an edge $e_{04}$ to the
complex $X$.  We denote the new complex by $X'$.  The set of data for
the new triangulation is
$$\Tet'(Y)=\Tet(Y) - \{S(0123),S(1234)\} \cup\{S(0234), S(0134),
S(0124)\},$$
$$\Edge'(Y)=\Edge(Y)\cup \{ e_{04} \}, \ \ \ \Ver'(Y)=\Ver(Y).$$  
Any flat labelling $l$ of $\Edge(T)$ induces a flat labelling $l'$ of
$\Edge(T')$.  Since any loop in $Y$ can be deformed so as not to
intersect $S(0123) \cup S(1234)$, $[\rho]$ is the same for $l$ and
$l'$.  Let $l_{\on{new}}^{(0)}$ denote the function of the lengths
$l_1,\ldots,l_N$ given by the implicit function theorem so that if
$l_{\on{new}}^{(0)}$ is the length of the edge $(v_0v_4)$, and $l_j$
are other lengths, then $\omega_{\on{new}} = 0$.  Let $l_{\on{new}}$
denote the length of the edge $(v_0v_4)$, and
\begin{equation} \label{new}
\ti{l}_{\on{new}} = l_{\on{new}} - l_{\on{new}}^{(0)} .\end{equation}
Since
$$0 = \f{\partial \omega_{\on{new}}(l_{\on{new}}^{(0)}) }{ \partial l_j} = \f{\partial
\omega_{\on{new}}}{ \partial l_{\on{new}}} \f{ \partial l_{\on{new}}^{(0)}}{\partial l_j} +
\f{\partial \omega_{\on{new}}}{\partial l_j} $$
we have 
$$ d \ti{l}_{\on{new}} = dl_{\on{new}} + \sum_j \frac{ \partial
\omega_{\on{new}}/\partial l_j}{\partial \omega_{\on{new}}/ \partial l_{\on{new}}}
\d l_j .$$
It follows that 
\begin{equation} \label{factor}
\left(\begin{array}{c} d\omega_{\rm new}\\ d\omega_1
\\ \vdots\cr d\omega_N
\end{array}\right) =
\left(\begin{array}{cc} \partial \omega_{\rm new} / 
\partial l_{\rm new} &
0 \  \cdots \ 0\vspace{1mm}\\
\begin{array}{c}
\partial \omega_1 / \partial l_{\rm new}\\
 \vdots \\
  \partial \omega_N / \partial l_{\rm new}\end{array}
 &
 \hbox{\Huge{$H$}}\end{array}\right)
\left(\begin{array}{c}
 d\tilde l_{\rm new}\\ dl_1\\ \vdots\\ dl_N\end{array}\right).
\end{equation}

Hence
\begin{equation}\label{new A}
H_{\rm new} =
\left(\begin{array}{cc}
\partial \omega_{\rm new} / \partial l_{\rm new} &
 0 \ \cdots \ 0 \vspace{1mm}\\
\begin{array}{c} \partial \omega_1 / \partial l_{\rm new}\\
 \vdots \\   \partial \omega_N / \partial l_{\rm new} \end{array}
& \hbox{\Huge{$H$}}\end{array}\right)
\left(
\begin{array}{cccc}
1 &\frac{ \partial
\omega_{\on{new}}/\partial l_1}{\partial \omega_{\on{new}}/ 
\partial l_{\on{new}}}
&  \ldots & \frac{ \partial
\omega_{\on{new}}/\partial l_N}{\partial \omega_{\on{new}}/ 
\partial l_{\on{new}}} \\
           &  1   &        &      \\
           &      & \ddots &  \hbox{\Huge 0}\\
           & \hbox{\Huge 0}  &        &   1 \end{array}\right).
\end{equation}

Since both matrices are block triangular,
$$ \det(H_{\on{new}}) = \frac{\partial \omega_{\on{new}}}{\partial l_{\on{new}}}
\det(H) .$$
The new triangulation has $\cC'=\cC \cup \{(v_0v_4)\}$, so that
$\ol{\cC'}=\Edge'-\cC'=\ol{\cC} .$ Invariance now follows from
\eqref{sjac}.

Next we show that $I(Y,[\rho])$ is independent of the choice
of $\cC$.  We write
$$l = (l',l''), \ \ \ \omega = (\omega', \omega'')$$ 
where $l'$ is the vector of edge lengths in $\cC$, and $l''$ the the
remaining edge lengths, and similarly for $\omega$.  Generically the
length $l''$ may be written as a function of $l'$, by requiring that
the defect angles $\omega = 0$.  With respect to this decomposition,
the matrix $H$ may be written in block diagonal form as follows.  Let
\begin{equation} \label{Ddef}
 D = \frac{\partial l_i''}{\partial l_j'}\end{equation}
denote the matrix of partial derivatives.  Define
$$\d \ti{l}' = \d l' + \frac{\partial l''}{\partial l'} \d l'' $$ 
similar to \eqref{new}.  It follows from the definition that
$$ \begin{pmatrix} \d \omega' \\ \d \omega'' \end{pmatrix}
= \begin{pmatrix} B & 0 \\ C & 0 \end{pmatrix}
\begin{pmatrix} \d \ti{l}' \\ \d l'' \end{pmatrix} $$
for some matrices $B,C$.  We have an equation similar to
\eqref{factor}
$$ H = \begin{pmatrix} B & 0 \\
                C & 0 \end{pmatrix} 
        \begin{pmatrix}I &  D \\
                0 & I  \end{pmatrix} = \begin{pmatrix} B & B D \\
                                       C & C D \end{pmatrix}  .$$
It follows from the fact that $H$ is symmetric that 
\begin{equation} \label{blockA} H =  \begin{pmatrix}   
H_{\cC \cC} 
&  H_{\cC \cC} D \\
D^T H_{\cC \cC} & D^T H_{\cC \cC} D \end{pmatrix}.
\end{equation}
Let $\cC'$ be a different maximal subset of edges, such that $
H_{\cC' \cC'}$ is non-degenerate.  Take $X \subset \cC$, $Y \subset
\ol{\cC}$ such that $|X|=|Y|$.  Set $\cC'=(\cC-X)\cup Y$.  From
\eqref{blockA} we see that
\begin{equation}
H_{\cC' \cC'}=\begin{pmatrix} H_{\cC-X \cC-X} & H_{\cC-X X}
  D_{XY} \\
D_{XY}^T H_{X \cC-X} & D_{XY}^T H_{XX} D_{XY}
\end{pmatrix},
\end{equation}
since $\frac{\partial \omega_i}{\partial l_j}=\frac{\partial \omega_j}
{\partial l_i}$. 
Thus, 
$$ H_{\cC' \cC'} =  F^T (H_{\cC \cC}) F $$
where $F$ is the matrix in block diagonal form with respect to the
decomposition $\cC = (\cC \ssm X) \bigcup X $ for columns and
$\cC' = (\cC' \ssm Y) \bigcup Y$ for rows
$$ F = \begin{pmatrix} I & 0 \\ 0 & D_{XY} \end{pmatrix}.$$
It follows that 
$$ \det( H_{\cC' \cC'}) = \det(H_{\cC \cC}) \det(F)^2
= \det(H_{\cC \cC}) \det(D_{XY})^2 .$$
Together with \eqref{Ddef} this implies that the differential form
$$\frac{\bigwedge_{e \in \ol{\cC}}dl_e}{\sqrt{\det(H_{\cC \cC})}}$$ 
in \eqref{tvsc} is the same for $\cC$ and $\cC'$.

To prove invariance under a $1$-$4$ move we will use the following
lemma, whose proof is left to the reader:
\begin{lemma}\label{integral}
 The integral over the region 
$$\{l_b,l_c : l_a\leq l_b+l_c, l_b \leq
l_a+l_c, l_c\leq l_a+l_b, l_a+l_b+l_c\leq 2\pi \}$$
\begin{equation} \label{delinfty}
\frac{1}{\sin(l_a)}\int\int\sin(l_b)\sin(l_c) \d l_b \d l_c =2 
\end{equation}
for any $l_a \in [0,\pi]$.
\end{lemma}

Let $S_{0123}$ be a tetrahedron with vertices $v_0, \ldots, v_3$ in
$Y$.  We consider the effect of adding an extra vertex $v_4$ in the
interior and replacing the tetrahedron $S_{0123}$ with the four
tetrahedra $ S_{1234}$, $S_{0234}$, $S_{0134}$, $S_{0124}$. We use the
notation $\tau_i$ for $S_{0 \ldots \hat{i} \ldots 4}$.  We have
$$\Ver'=\Ver \cup \{ v_4 \}, \ \ \Tet'=(\Tet -
\{S(0123)\}) \cup \{S(1234),S(0234),S(0134),S(0124) \} $$
$$\Edge'=\Edge \cup \{
e_{04},e_{14},e_{24},e_{34}\}.
$$
Also $\cC' = \cC \cup e_{34}$ since adding any other edge would allow
a deformation of the new vertex changing only the lengths of edges in
$\cC'$.  Hence
$$\ol{\cC'}=\ol{\cC}\cup\{e_{04},e_{14},e_{24}\}.$$ 
Exactly the same argument as in the 2-3 case shows that
$$\det (H_{\cC' \cC'})=\frac{\partial \omega_{34}}{\partial
l_{34}}\det ( H_{\cC \cC}).$$
Hence
$$I(Y')=(\frac{1}{2\pi})^{(\# \Ver +1)}\int_{\mL'}
\frac{\prod_{\tau \in \Tet-\{S(0123)\}}(G(l(\tau))^{-1/4})
\prod_{e \in \Edge'}\sin( l_e)\bigwedge_{e \in
\ol{\cC'}}dl_e}
{(G_0G_1G_2G_3)^{1/4}\sqrt{(\det H
 _{\cC \cC})(\frac{\partial \omega_{34}}{\partial l_{34}})}}.$$
After substituting the Jacobian 
$$\frac{\partial \omega_{34}}{\partial l_{34}}=\sin^2(l_{34})
\sqrt{\frac{G_3G_4}{G_0G_1G_2}}$$ 
we need to compute the integral
$$\int \frac{\sin( l_{04})\sin( l_{14})\sin(
l_{24})dl_{04}dl_{14}dl_{24}}
{\sqrt{G_3}}.$$
The equations \eqref{normalization} and \eqref{delinfty} give
$$\frac{\pi}{\sin( l_{12})}\int \sin( l_{14})\sin(
l_{24})dl_{14}dl_{24}=2\pi.$$ 
This cancels with the extra factor of $2\pi$ in the coefficient, and
completes the proof that $I(Y,[\rho])$ is invariant under the Pachner
moves, assuming it converges.

\subsection{An acyclic complex and its torsion}

In this section we relate the determinant appearing in $I(Y,[\rho])$
to the torsion of an acyclic complex, following Korepanov
\cite{ik:cp}.  Recall the infinitesimal action of the group of gauge
transformations $\Map(Y,SO(4))$ on the space of connections
$\Omega^1(Y,\so(4))$ at a connection $A$ is given by
\begin{equation} \label{infin} \Omega^0(Y,\so(4)) \to \Omega^1(Y,\so(4)),
  \ \ \xi \mapsto - d_A \xi \end{equation}
where $d_A$ is the associated covariant derivative.  Hence the
infinitesimal stabilizer of $A$ is

$$ \Omega^0(Y,\so(4))_A = H_0(d_A) .$$

Let $\SO(4)_\rho$ denote the stabilizer of $\rho: \pi_1(Y) \to SO(4)$,
and $\so(4)_\rho$ its Lie algebra.  If $A$ is a flat connection
defining the holonomy representation $\rho$, then evaluation at the
identity induces an isomorphism
$$ \Map(Y,SO(4))_A \to SO(4)_\rho .$$
Hence $ H^0(d_A)$ is isomorphic to $\so(4)_\rho$.  Let 
$$ h^0(d_A) = \dim(\so(4)_\rho) = \dim(H^0(d_A)) .$$ 
The cohomology group $H^1(d_A)$ parameterizes first-order deformations
of $\rho$; in particular, if $H^1(d_A) = 0$ then $[\rho]$ is isolated
in $R(Y,SO(4))$.  Suppose that $H^1(d_A) = 0$.  Let
$$ V = \Map(\Ver(\ti{Y}),S^3)^{\pi_1(Y)} $$
denote the space of maps invariant under $\pi_1(Y)$, acting on
$\Ver(\ti{Y})$ by deck transformations and $S^3$ via the
representation $\rho$.  Let
$$ E = \Map(\Edge(\ti{Y}),[0,\pi])^{\pi_1(Y)} = \Map(\Edge(Y), [0,\pi]) $$
and $ \delta: V \to E $
the map taking edge lengths of edges.  Let
$ \omega: E \to E $
be the map which assigns to a set of edge lengths the set of defect
angles.  The action of $SO(4)_\rho$ on $V$ induces a map
$$ \lambda: \ \so(4)_\rho \to \Vect(V) .$$
Evaluating the vector field at $p \in V$ gives
$$ \lambda_p: \ \so(4)_\rho \to T_p V .$$
For any $p \in V$, let $l = \delta(p), \ \hat{l} = \omega(l) $ and
$\hat{p}$ any point in $\delta^{-1}(\hat{l})$.  Consider the sequence
\begin{equation}\label{complex}
0 \to \so(4)_\rho \to T_p V \to T_l E \to T_{\hat{l}} E \to T
_{\hat{p}} V \to \so(4)_\rho \to 0
\end{equation}
with maps $\lambda_p,D_p \delta,\Hess,D_p
\delta^T,\lambda_{\hat{p}}^T$.  It follows from the fact that $\Hess$
is symmetric and a straight-forward calculation that the sequence
\eqref{complex} is exact, that is, \eqref{complex} is an acyclic
complex.  Let $\tau(l,s)$ denote the torsion, which is defined as
follows.  Let $\Ver'(Y)$ denote a maximal subset of the space of
vertices so that $\delta$ is injective on the corresponding subspace
of $T_p V$. Let $\delta'$ denote the restriction of $\delta$ to
$\Ver'(Y)$, followed by projection onto the subspace of $T_l E$
corresponding to the complement of $\cC$.  Then
$$ \tau(l,s) =\det(\lambda)^{-2} \det(\delta')^2 \det(\Hess_{\cC \cC})^{-1}.$$%

\section{Computations of the Invariant for the sphere $S^3$}\label{comp}

A triangulation of $S^3$ consists of the following data:
$$\Ver=\{ 0,1,2,3,4\} ,$$
$$\Edge=\{ 01,02,03,04,12,13,14,23,24,34\} ,$$
$$\Face=\{ 012,023,013,124,123,134,234,014,024,034\} ,$$ 
$$\Tet=\{ 0123,1234,0124,0234,0134 \} .$$
Since $S^3$ is simply-connected, the representation variety $R(S^3,
SO(4))$ is trivial.  So, $I(S^3)= I(S^3,[1]).$
Using the acyclic complex in the previous section, we find that the
rank of $\cC$ is 1.  Since there is no distinction among the edges, we
choose $\cC=\{04\}$.  Thus, $\ol{\cC}=\Edge - \{04\}.$ Denote by
$G_{i}$ the determinant of the Gram matrix associated with the
tetrahedron $(0 \ldots \hat{i} \ldots 4)$.  We must compute
\begin{equation}\label{invariant}
(\frac{1}{2\pi})^5
\int_{\mL_{\flat, \on{I}}}\frac{\prod_{e \in \Edge}\sin(
l_e)\bigwedge_{e \in
\ol{\cC}}dl_e}{(G_0G_1G_2G_3G_4)^{1/4}\sqrt{\det{H_{_\cC \cC}}}}.
\end{equation}
Note that around the edge $(04)$, there are three tetrahedra
$(0234),(0124),(0134)$.  When these tetrahedra match with each other
in $S^3$ under the curvature zero condition around the edge $(04)$, we
have two tetrahedra $(0123)$ and $(1234)$ as well.  In other words, we
are in the situation where the spherical Jacobian \eqref{sjac} is
equal to $H_{\cC \cC}$.  So, the integral reduces to
$$(\frac{1}{2\pi})^5 \int \frac{\prod_{e \in \ol{\cC} }\sin(
l_e)\bigwedge_{e \in \ol{\cC} }dl_e}{\sqrt{G_0G_4}}.$$
Apply the orthogonality identity \eqref{normalization} to the
tetrahedra $(1234)$ and $(0123)$ respectively and integrate the rest
from $0$ to $\pi$ in each variable.  Then, the integral
\eqref{invariant} is computed to be $\frac{1}{\pi^3}.$ There are $2^5$
ways of assigning signs to each tetrahedron in the triangulation, but
the above argument is applied to each assignment of the sign.
Therefore,
$$I(S^3 )=\frac{2^5}{\pi^3}.$$

\section{Remarks on the semiclassical limit of Turaev-Viro}

Throughout this section we assume that $Y$ is a rational homology
sphere.  The stationary phase approximation to the Chern-Simons path
integral predicts \cite{fg:91}
$$ Z(Y,r) \sim \hh r^{-\hh h^0(d_A)}
e^{-3\pi i /4} \sum_{[A] \in R(Y,SU(2))}
\sqrt{\tau(A)} e^{-2\pi i I_A/4} e^{2\pi i CS(A,r)} $$
where $\tau(A)$ is the torsion of $A$, $I_A$ is the spectral
flow, and $CS(A,r)$ the Chern-Simons invariant at level $r$
$$ CS(A,r) = \frac{r}{8 \pi^2} \int_Y \Tr(A \wedge \d A + \frac{2}{3}
A \wedge A \wedge A \} .$$
We write any $SO(4)$ connection as a pair of $SU(2)$-connections.  The
norm-square of the asymptotic formula for $Z(Y,r)$ is
\begin{equation} \label{tvasym}
 TV(Y,r) \sim \qq \sum_{[A] \in R(Y,SO(4))} r^{- \hh
 h^0(d_A)} \sqrt{\tau(A_1) \tau(A_2)} e^{-2\pi (I_{A_1} - I_{A_2})/4}
 e^{2\pi i (CS(A_1,r) - CS(A_2,r))} \end{equation}
where $A = (A_1,A_2)$.  

\subsection{The leading power of $r$}

It follows from $\Delta(r) = r \sin( \pi/r)^{-2} $ that $\Delta \sim
\frac{r^3}{\pi^2}$ as $r \to \infty$.  Let $t,e,v$ denote the size of
the sets $\Tet(Y)$, $\Edge(Y)$, $\Ver(Y)$.  Collecting together the
powers of $r$ in the asymptotic $6j$ formula \eqref{S^3-asymp}, the
definition of the Turaev-Viro invariant \eqref{tv}, and the acyclicity
of \eqref{complex} we obtain the prediction for leading power of $r$
in the Turaev-Viro invariant
$$ - \thh v + \thh e - \thh t - \hh h^0(d_A) = -
  \hh h^0(d_A)  .$$
This agrees with the prediction in \eqref{tvasym}.

\subsection{The Volumes/Chern-Simons invariant}

The terms $\exp(\pm i \phi)$ appearing in the stationary phase
approximation to Turaev-Viro lead to a factor 
$$ \exp \left( \frac{i}{\pi} \sum_{\tau \in \Tet(Y)} \pm \Vol(\tau) \right) .$$
Let $\phi:\tilde{Y} \to S^3$ denote the developing map as in
Proposition \ref{varphi}.  Let $\d \Vol(S^3)$ denote the volume form
on $S^3$ so that $\int_{S^3}\d \Vol(S^3)= 2 \pi^2$.

Let $\pi: \ SO(4) \to S^3$ denote the map given by action on
$(1,0,0)$.  We have $ \pi^* \d \Vol(S^3) = 2 \pi^2 \chi $ where
$\chi = (\alpha,[\alpha,\alpha]) \in \Omega^3(SO(4))$ is the
Chern-Simons three-form on $SO(4)$ with $\alpha \in
\Omega^1(SO(4),\so(4))$ the left Maurer-Cartan form and $(\ , \ )$ the
inner product equal to the basic inner product on one $\su(2)$-factor
and minus the basic inner product on the other.  Let $A = (A_1,A_2)$
be an $SU(2)^2$ connection on $Y$ with holonomy representation $\rho$
and $g: \ti{Y} \to SU(2)^2$ a gauge transformation trivializing the
lift $\ti{A}$ of $A$ to $\ti{Y}$.  For any $\gamma \in \pi_1(Y)$, we
have $\gamma^* g = \rho(\gamma) g$.  This implies that $g^{-1} \cdot
\phi$ is $\pi_1$-invariant, and hence descends to a map $Y \to S^3$.
Hence
\begin{eqnarray*}
 \frac{1}{\pi} \sum_{\tau \in \Tet(Y)} \pm \Vol(\tau) &=& 
\frac{1}{\pi}   (\# \pi_1(Y))^{-1}
 \int_{\ti{Y}} \phi^* \d \Vol(S^3) \\
&=& \frac{1}{\pi}   (\# \pi_1(Y))^{-1} \int_{\ti{Y}}
 g^* \pi^* \d \Vol(S^3)   \ \ \ \ \text{mod} \ 2\pi\Z \\
&=&  2 \pi ( \# \pi_1(Y))^{-1} \int_{\ti{Y}} g^*
 \chi   \ \ \ \ \text{mod} \ 2\pi\Z \\
&=&   2 \pi (\# \pi_1(Y))^{-1} ( CS(\ti{A}_1) - CS(\ti{A}_2))  \ \ \ \
 \text{mod} \ 2\pi\Z \\
&=&  2 \pi ( CS(A_1) - CS(A_2))  \ \ \ \ \text{mod} \ 2\pi\Z \end{eqnarray*}
which also matches \eqref{tvasym}.

\subsection{The Maslov indices and torsion}

Each tetrahedron contributes $\exp( \pm \pi i /4)$ from the formula
\eqref{S^3-asymp}.  Stationary phase leads to a factor $\exp( \pi i
\on{sign}(H_{\cC \cC})/4)$.  It seems natural to conjecture that these
combine to the spectral flow factor $\exp( 2 \pi i I_A/4)$ in the
Freed-Gompf formula.  One expects the torsion to correspond to our
three-manifold invariant.  However, it is not clear to us how to
perform the integral over flat labellings.

Revised June 12, 2004.



\end{document}